\documentclass[12pt]{article}

\usepackage{amsmath, epsfig, cite}
\usepackage{amssymb}
\usepackage{amsfonts}
\usepackage{latexsym}
\usepackage{float}
\usepackage{color}
\usepackage{mathrsfs}
\usepackage{booktabs}
\usepackage{graphicx}
\usepackage{caption}

\newtheorem{thm}{Theorem}[section]

\numberwithin{equation}{section}

\newcommand{\qed}{{\hfill$\square$}\medskip}
\UseRawInputEncoding
\begin{document}

\begin{center}
{\large\bf A combinatorial approach to Berkovich type identities}
\end{center}

\vskip 2mm \centerline{Ji-Cai Liu}
\begin{center}
{\footnotesize Department of Mathematics, Wenzhou University, Wenzhou 325035, PR China\\
{\tt jcliu2016@gmail.com} }
\end{center}


\vskip 0.7cm \noindent{\bf Abstract.}
Motivated by Berkovich's nine $q$-binomial identities involving the Legendre symbol $(\frac{d}{3})$,
we establish a unified form of $q$-binomial identities of this type through a combinatorial approach.
This unified form includes Berkovich's nine identities as special cases. Many such identities can be also deduced from this unified form.

\vskip 2mm
\noindent{\it MR Subject Classifications}: 05A19, 05A10, 05A17

\section{Introduction}
Berkovich and Uncu \cite{bu-jmaa-2022} established many interesting $q$-binomial identities involving the Legendre symbol $(\frac{d}{3})$, which is given by
\begin{align*}
\left(\frac{d}{3}\right)=
\begin{cases}
0\quad&\text{if $d\equiv 0\pmod{3}$},\\
1\quad&\text{if $d\equiv 1\pmod{3}$},\\
-1\quad&\text{if $d\equiv -1\pmod{3}$}.
\end{cases}
\end{align*}
For instance, they proved that for non-negative integers $L$,
\begin{align*}
\sum_{j=-L}^L\left(\frac{j+1}{3}\right)q^{j^2}{2L\brack L+j}=
\sum_{m,n\ge 0}\frac{q^{2m^2+6mn+6n^2}(q;q)_L}{(q;q)_m(q^3;q^3)_n(q;q)_{L-3n-2m}}.
\end{align*}
Here and throughout the paper, the $q$-shifted factorial is defined by $(a;q)_n=(1-a)(1-aq)\cdots (1-aq^{n-1})$ for $n\ge 1$ and $(a;q)_0=1$. The $q$-binomial coefficients are defined as
\begin{align*}
{n\brack k}
=\begin{cases}
\displaystyle\frac{(q;q)_n}{(q;q)_k(q;q)_{n-k}} &\text{if $0\leqslant k\leqslant n$},\\[10pt]
0 &\text{otherwise.}
\end{cases}
\end{align*}

Recently, Berkovich \cite{berkovich-dm-2024} further investigated more $q$-binomial identities involving the Legendre symbol $(\frac{d}{3})$. Six amazing identities are established by Berkovich \cite[Theorems 2.1--2.6]{berkovich-dm-2024}, which are listed as follows (equivalent forms): For positive integers $L$,

\begin{align}
&\sum_{j=-L}^{L}\left(\frac{j}{3}\right)q^{\frac{j^2-3j}{2}}{2L\brack L+j}=
\frac{(-1;q^3)_{L-1}}{(-1;q)_{L-1}}q^{L-2}(1-q^L)(1+q+q^2),\label{berk-1}\\[5pt]
&\sum_{j=-L}^{L+1}\left(\frac{j}{3}\right)q^{\frac{j^2-3j}{2}}{2L+1\brack L+j}=
\frac{(-1;q^3)_{L}}{(-1;q)_{L}}q^{L-1}(1+q-q^{L+1}),\label{berk-2}\\[5pt]
&\sum_{j=-L-1}^{L}(-1)^j\left(\frac{j+2}{3}\right)q^{\frac{j^2}{2}}{2L+1\brack L+1+j}=
-\frac{(q^{3/2};q^3)_L}{(q^{1/2};q)_{L+1}}(1-q^{2L+1}),\label{berk-3}\\[5pt]
&\sum_{j=-L}^{L}(-1)^j\left(\frac{j+1}{3}\right)q^{\frac{j^2}{2}}{2L\brack L+j}=
\frac{(q^{3/2};q^3)_L}{(q^{1/2};q)_{L}},\label{berk-4}\\[5pt]
&\sum_{j=-L-1}^{L}(-1)^j\left(\frac{j+1}{3}\right)q^{\frac{j^2}{2}}{2L+1\brack L+1+j}=
\frac{(q^{3/2};q^3)_L}{(q^{1/2};q)_{L}},\label{berk-5}\\[5pt]
&\sum_{j=-L-1}^{L}(-1)^j\left(\frac{j}{3}\right)q^{\frac{j^2}{2}}{2L+1\brack L+1+j}=
\frac{(q^{3/2};q^3)_L}{(q^{1/2};q)_{L}}q^{L+1/2}.\label{berk-6}
\end{align}

Berkovich \cite{berkovich-dm-2024} also listed three additional $q$-binomial identities of the same type without proof:
\begin{align}
&\sum_{j=-L-1}^{L}(-1)^j\left(\frac{j+1}{3}\right)q^{\frac{j^2-j}{2}}{2L+1\brack L+1+j}=
\frac{(q^3;q^3)_{L-1}}{(q;q)_{L-1}}(2+q^L+q^{L+1}-q^{2L+1}),\label{berk-7}\\[5pt]
&\sum_{j=-L-1}^{L}(-1)^j\left(\frac{j}{3}\right)q^{\frac{j^2-j}{2}}{2L+1\brack L+1+j}=
\frac{(q^3;q^3)_{L-1}}{(q;q)_{L-1}}(-1+q^L+q^{L+1}+2q^{2L+1}),\label{berk-8}\\[5pt]
&\sum_{j=-L-1}^{L}(-1)^{j+1}\left(\frac{j+2}{3}\right)q^{\frac{j^2-j}{2}}{2L+1\brack L+1+j}=
\frac{(q^3;q^3)_{L-1}}{(q;q)_{L-1}}(1+2q^L+2q^{L+1}+q^{2L+1}).\label{berk-9}
\end{align}

The important ingredients in Berkovich's proof include the polynomial analogue of the Jacobi triple product identity \cite[page 49]{andrews-b-1998} and the following Eisenstein formula for the Legendre symbol:
\begin{align*}
\left(\frac{j}{p}\right)=\prod_{n=1}^{(p-1)/2}\frac{\sin(2nj\pi/p)}{\sin(2n\pi/p)},
\end{align*}
for an odd prime $p$.

Note that all of the nine identities due to Berkovich \cite{berkovich-dm-2024} are in the unified form:
\begin{align*}
\sum_{k=-n}^m \varepsilon^{k+d} \left(\frac{k+d}{3}\right) q^{\frac{k^2+sk}{2}}{n+m\brack n+k}=
\text{$q$-product}\times \text{sum of finite terms in $q$},
\end{align*}
where $\varepsilon\in \{1,-1\},d\in\{0,1,2\}$ and $s\in \mathbb{Z}$.

The motivation of the paper is to establish the above unified form through a combinatorial approach.

For a formal Laurent series $f(z)=\sum_{k=-\infty}^{\infty} c_kz^k$, the operator $\mathcal{L}_z$ with respect to the Legendre symbol modulo $3$ is defined by
\begin{align*}
\mathcal{L}_z\left(f(z)\right)=\sum_{k=-\infty}^{\infty} \left(\frac{k}{3}\right)c_k.
\end{align*}
It is clear that
\begin{align*}
\mathcal{L}_z\left(f(z)\right)=\sum_{k=-\infty}^{\infty}c_{3k+1}-\sum_{k=-\infty}^{\infty}c_{3k+2}.
\end{align*}

The main results of the paper consist of the following two theorem.
\begin{thm}[Odd case]\label{t-1}
Let $\varepsilon\in \{1,-1\},d\in\{0,1,2\}$ and $s<0$ be an odd integer. For positive integers $n,m$ with $(1-s)/2\le m\le n-s$, we have
\begin{align}
&\sum_{k=-n}^m \varepsilon^{k+d} \left(\frac{k+d}{3}\right) q^{\frac{k^2+sk}{2}}{n+m\brack n+k}\notag\\[5pt]
&=\frac{(-\varepsilon q^3;q^3)_{m+(s-1)/2}}{(-\varepsilon q;q)_{m+(s-1)/2}}\mathcal{L}_z\left(q^{\frac{1-s^2}{8}}(\varepsilon z)^{d-\frac{s+1}{2}}(1+\varepsilon z)(-\varepsilon q^{m+\frac{s+1}{2}}/z;q)_{n-m-s}\right).\label{a-1}
\end{align}
\end{thm}

\begin{thm}[Even case]\label{t-2}
Let $\varepsilon\in \{1,-1\},d\in\{0,1,2\}$ and $s\le 0$ be an even integer. For positive integers $n,m$ with $-s/2\le m\le n-s$, we have
\begin{align}
&\sum_{k=-n}^m \varepsilon^{k+d} \left(\frac{k+d}{3}\right) q^{\frac{k^2+sk}{2}}{n+m\brack n+k}\notag\\[5pt]
&=\frac{(-\varepsilon q^{\frac{3}{2}};q^3)_{m+s/2}}{(-\varepsilon q^{\frac{1}{2}};q)_{m+s/2}}\mathcal{L}_z\left(q^{-\frac{s^2}{8}} (\varepsilon z)^{d-\frac{s}{2}}(-\varepsilon q^{m+\frac{s+1}{2}}/z;q)_{n-m-s}\right).\label{a-2}
\end{align}
\end{thm}

The rest of the paper is organized as follows. In the next section, we list some special examples of Theorems \ref{t-1} and \ref{t-2} which include Berkovich's nine identities and seven new simple examples.
In Section 3, we establish two combinatorial ingredients which will be utilized in the proof of Theorems \ref{t-1} and \ref{t-2}. We prove Theorems \ref{t-1} and \ref{t-2} in the last section.

\section{Special examples}
We write \eqref{a-1} and \eqref{a-2} in a unified form:
\begin{align*}
\sum_{k=-n}^m \varepsilon^{k+d} \left(\frac{k+d}{3}\right) q^{\frac{k^2+sk}{2}}{n+m\brack n+k}
=R_{\varepsilon,s,m}(q)\mathcal{L}_z(*),
\end{align*}
where $R_{\varepsilon,s,m}(q)$ denotes the $q$-product part and $\mathcal{L}_z(*)$ denotes
the computed result by the operator $\mathcal{L}_z$. If $s$ and $n-m$ are fixed integers, one can
easily calculate $\mathcal{L}_z(*)$.

Berkovich's nine identities \eqref{berk-1}--\eqref{berk-9} are included in Theorems \ref{t-1} and \ref{t-2} as special cases. The values $\mathcal{L}_z(*)$ for the nine identities due to Berkovich \cite{berkovich-dm-2024} are listed in the following table.
\begin{table}[H]
\caption*{\text{Nine known examples of Theorems \ref{t-1} and \ref{t-2}}}
\centering
\begin{tabular}{cccc}
\toprule
$(\varepsilon,d,s,n,m)$&$\mathcal{L}_z(*)$&$R_{\varepsilon,s,m}(q)$&Identities\\
\midrule
$(1,0,-3,L,L)$&$q^{L-2}(1-q^L)(1+q+q^2)$&$\frac{(-q^3;q^3)_{L-2}}{(-q;q)_{L-2}}$& \eqref{berk-1}\\[10pt]
$(1,0,-3,L,L+1)$&$q^{L-1}(1+q-q^{L+1})$&$\frac{(-q^3;q^3)_{L-1}}{(-q;q)_{L-1}}$& \eqref{berk-2}\\[10pt]
$(-1,2,0,L+1,L)$&$-1-q^{L+1/2}$&$\frac{( q^{3/2};q^3)_{L}}{( q^{1/2};q)_{L}}$&\eqref{berk-3}\\[10pt]
$(-1,1,0,L,L)$&$-1$&$\frac{( q^{3/2};q^3)_{L}}{( q^{1/2};q)_{L}}$&\eqref{berk-4}\\[10pt]
$(-1,1,0,L+1,L)$&$-1$&$\frac{( q^{3/2};q^3)_{L}}{( q^{1/2};q)_{L}}$&\eqref{berk-5}\\[10pt]
$(-1,0,0,L+1,L)$&$q^{L+1/2}$&$\frac{( q^{3/2};q^3)_{L}}{( q^{1/2};q)_{L}}$&\eqref{berk-6}\\[10pt]
$(-1,1,-1,L+1,L)$&$q^{2L+1}-q^{L+1}-q^L-2$&$\frac{(q^3;q^3)_{L-1}}{(q;q)_{L-1}}$&\eqref{berk-7}\\[10pt]
$(-1,0,-1,L+1,L)$&$2q^{2L+1}+q^{L+1}+q^L-1$&$\frac{(q^3;q^3)_{L-1}}{(q;q)_{L-1}}$&\eqref{berk-8}\\[10pt]
$(-1,2,-1,L+1,L)$&$-q^{2L+1}-2q^{L+1}-2q^L-1$&$\frac{(q^3;q^3)_{L-1}}{(q;q)_{L-1}}$&\eqref{berk-9}\\
\bottomrule
\end{tabular}
\end{table}

For $-5\le s\le 0$ and $0\le n-m\le 5$, we search for quintuples $(\varepsilon,d,s,n,m)$
such that $\mathcal{L}_z(*)$ is a constant. Seven new simple examples are found, the values $\mathcal{L}_z(*)$ of which are constants.

\begin{table}[H]
\caption*{\text{Seven new simple examples}}
\centering
\begin{tabular}{ccc}
\toprule
$(\varepsilon,d,s,n,m)$&$\mathcal{L}_z(*)$&$R_{\varepsilon,s,m}(q)$\\
\midrule
$(1, 0, 0, L, L)$&$0$&$\frac{(-q^{3/2};q^3)_{L}}{(-q^{1/2};q)_{L}}$ \\[10pt]
$(1, 1, 0, L, L)$&$1$&$\frac{(-q^{3/2};q^3)_{L}}{(-q^{1/2};q)_{L}}$ \\[10pt]
$(1, 2, 0, L, L)$&$-1$&$\frac{(-q^{3/2};q^3)_{L}}{(-q^{1/2};q)_{L}}$ \\[10pt]
$(-1, 0, 0, L, L)$&$0$&$\frac{(q^{3/2};q^3)_{L}}{(q^{1/2};q)_{L}}$\\[10pt]
$(-1, 2, 0, L, L)$&$-1$&$\frac{(q^{3/2};q^3)_{L}}{(q^{1/2};q)_{L}}$\\[10pt]
$(1, 2, -1, L, L)$&$-1$&$\frac{(-q^3;q^3)_{L-1}}{(-q;q)_{L-1}}$\\[10pt]
$(1, 1, 0, L+1, L)$&$1$&$\frac{(-q^{3/2};q^3)_{L}}{(-q^{1/2};q)_{L}}$\\
\bottomrule
\end{tabular}
\end{table}

\section{Two combinatorial ingredients}
In this section, we will establish two important combinatorial ingredients in the proof of Theorems \ref{t-1} and \ref{t-2}, which are interesting by themselves.
\begin{thm}
Let $n,m$ be positive integers and $s$ be an integer. Then
\begin{align}
\sum_{k=-n}^m z^{n+k} q^{\frac{k^2+sk}{2}}{n+m\brack n+k}
=q^{\frac{n(n-s)}{2}}(-zq^{-n+\frac{s+1}{2}};q)_{n+m}.\label{c-new-1}
\end{align}
\end{thm}

{\noindent \it Proof.}
Recall that a partition of a positive integer $n$ is a finite nonincreasing sequence of positive integers $\lambda_1,\lambda_2,\cdots,\lambda_r$ such that $\sum_{i=1}^r\lambda_i=n$.
By \cite[Theorem 3.1, page 33]{andrews-b-1998}, we have
\begin{align*}
{n\brack k}=\sum_{N\ge 0}p(n-k,k,N)q^N,
\end{align*}
where $p(k,n-k,N)$ is the number of partition of $N$ into at most $k$ parts and each part $\le n-k$.
Since ${k+1\choose 2}=1+2+\cdots+k$, we have
\begin{align}
{n\brack k}q^{{k+1\choose 2}}=\sum_{N\ge 0}p_{d}(k,n,N)q^N,\label{c-1}
\end{align}
where $p_{d}(k,n,N)$ is the number of partition of $N$ into $k$ distinct parts and each part $\le n$.

Let $[n]=\{1,\cdots,n\}$ and $\mathscr{G}=\{A\subseteq [n]\}$. For any $A\in \mathscr{G}$, we associate $A$ with a weight $||A||=\sum_{a\in A}a$. We rewrite \eqref{c-1} as
\begin{align}
{n\brack k}q^{{k+1\choose 2}}=\sum_{A\subseteq [n]\atop \#A=k}q^{||A||}.\label{c-2}
\end{align}

Let $S=\{-n+(s+1)/2,\cdots,m+(s-1)/2\}$, which is obtained by $[n+m]$ by a shift $-n+(s-1)/2$, and $\mathscr{F}=\{A\subseteq S\}$. By \eqref{c-2}, we have
\begin{align}
\sum_{A\in \mathscr{F}}z^{\#A}q^{||A||}&=\sum_{k=-n}^{m}z^{n+k}\sum_{A\subseteq S\atop \#A=n+k}q^{||A||}\notag\\[5pt]
&=\sum_{k=-n}^{m}z^{n+k} {n+m\brack n+k}q^{{n+k+1\choose 2}+(n+k)(-n+\frac{s-1}{2})}\notag\\[5pt]
&=q^{-\frac{n(n-s)}{2}}\sum_{k=-n}^{m}z^{n+k} {n+m\brack n+k}q^{\frac{k^2+sk}{2}}.\label{c-3}
\end{align}

On the other hand, we have
\begin{align}
\sum_{A\in \mathscr{F}}z^{\#A}q^{||A||}
&=(1+zq^{-n+\frac{s+1}{2}})\cdots (1+zq^{m+\frac{s-1}{2}})\notag\\[5pt]
&=(-zq^{-n+\frac{s+1}{2}};q)_{n+m}.\label{c-4}
\end{align}

Finally, combining \eqref{c-3} and \eqref{c-4}, we arrive at \eqref{c-new-1}.
\qed

\begin{thm}
For $\varepsilon\in \{1,-1\}$, positive integers $n$ and integers $d$, we have
\begin{align}
\mathcal{L}_z\left(z^d (\varepsilon zq;q)_n (\varepsilon q/z;q)_n\right)=\left(\frac{d}{3}\right)\frac{(\varepsilon q^3;q^3)_n}{(\varepsilon q;q)_n},\label{c-6}
\end{align}
and
\begin{align}
\mathcal{L}_z\left(z^d (\varepsilon zq^{\frac{1}{2}};q)_n (\varepsilon q^{\frac{1}{2}}/z;q)_n\right)=\left(\frac{d}{3}\right)\frac{(\varepsilon q^{\frac{3}{2}};q^3)_n}{(\varepsilon q^{\frac{1}{2}};q)_n}.\label{c-7}
\end{align}
\end{thm}

{\noindent \it Proof.}
Since $(\frac{d+3l}{3})=(\frac{d}{3})$, it suffices to show that \eqref{c-6} and \eqref{c-7} are true for
$d\in\{0,1,2\}$. Next, we shall prove the equivalent forms of \eqref{c-6} and \eqref{c-7}:
\begin{align}
\mathcal{L}_z\left(z^d (-\varepsilon zq;q)_n (-\varepsilon q/z;q)_n\right)=\left(\frac{d}{3}\right)\frac{(-\varepsilon q^3;q^3)_n}{(-\varepsilon q;q)_n},\label{c-6-new}
\end{align}
and
\begin{align}
\mathcal{L}_z\left(z^d (-\varepsilon zq^{\frac{1}{2}};q)_n (-\varepsilon q^{\frac{1}{2}}/z;q)_n\right)=\left(\frac{d}{3}\right)\frac{(-\varepsilon q^{\frac{3}{2}};q^3)_n}{(-\varepsilon q^{\frac{1}{2}};q)_n},\label{c-7-new}
\end{align}
for $d\in\{0,1,2\}$.

Let $S_n=\{\varepsilon zq,\cdots, \varepsilon zq^n,\varepsilon q/z,\cdots,\varepsilon q^n/z\}$ and
$\mathscr{H}_n=\{A\subseteq S_n\}$.
For any $A\in \mathscr{H}_n$, we associate $A$ with a weight and an index as follows:
\begin{align*}
&||A||_z=\text{the coefficient of $z$ in the product}\prod_{a\in A}a,\\[5pt]
&\text{index}_z(A)=\text{the index of $z$ in the product}\prod_{a\in A}a.
\end{align*}
For $r\in \{0,1,2\}$, let
\begin{align*}
\mathscr{H}_n^{(r)}=\{A\in \mathscr{H}_n|\quad\text{index}_z(A) \equiv r \pmod{3}\},
\end{align*}
and
\begin{align*}
C_n^{(r)}=\sum_{A\in \mathscr{H}_n^{(r)}} ||A||_z.
\end{align*}
We first proved that
\begin{align}
C_n^{(1)}=C_n^{(2)}.\label{c-8}
\end{align}

For $a=\varepsilon zq^k\in S_n$, let $a'=\varepsilon q^k/z$, and for $a=\varepsilon q^k/z\in S_n$, let $a'=\varepsilon zq^k$.
We define the involution $\sigma$ on $\mathscr{H}_n$ as follows:
\begin{align*}
\sigma(A)=\{a| a'\in A \}.
\end{align*}
For example, if $A=\{a,a',b,c\}$, then $\sigma(A)=\{a,a',b',c'\}$.
It is clear that for any $A\in \mathscr{H}_n$,
\begin{align}
&||\sigma(A)||_z=||A||_z,\label{c-9}\\[5pt]
&\text{index}_z(\sigma(A))=-\text{index}_z(A).\label{c-10}
\end{align}
It follows from \eqref{c-9} and \eqref{c-10} that $\sigma$ is a bijection between $\mathscr{H}_n^{(1)}$ and $\mathscr{H}_n^{(2)}$, and so $C_n^{(1)}=C_n^{(2)}$.

Note that
\begin{align*}
\mathscr{H}_n^{(0)}=&\mathscr{H}_{n-1}^{(0)}\biguplus \{A\cup \{\varepsilon q^n/z\}|A\in \mathscr{H}_{n-1}^{(1)}\} \biguplus \{A\cup \{\varepsilon zq^n\}|A\in \mathscr{H}_{n-1}^{(2)}\}\\[5pt]
&\biguplus \{A\cup \{\varepsilon zq^n,\varepsilon q^n/z\}|A\in \mathscr{H}_{n-1}^{(0)}\},\\[5pt]
\mathscr{H}_n^{(1)}=&\mathscr{H}_{n-1}^{(1)}\biguplus \{A\cup \{\varepsilon zq^n\}|A\in \mathscr{H}_{n-1}^{(0)}\} \biguplus \{A\cup \{\varepsilon q^n/z\}|A\in \mathscr{H}_{n-1}^{(2)}\}\\[5pt]
&\biguplus \{A\cup \{\varepsilon zq^n,\varepsilon q^n/z\}|A\in \mathscr{H}_{n-1}^{(1)}\}.
\end{align*}
It follows that
\begin{align}
&C_n^{(0)}=(1+(\varepsilon q^n)^2)C_{n-1}^{(0)}+\varepsilon q^n(C_{n-1}^{(1)}+C_{n-1}^{(2)}),\label{c-11}\\[5pt]
&C_n^{(1)}=(1+(\varepsilon q^n)^2)C_{n-1}^{(1)}+\varepsilon q^n(C_{n-1}^{(0)}+C_{n-1}^{(2)}).\label{c-12}
\end{align}
Combining \eqref{c-8}, \eqref{c-11} and \eqref{c-12}, we arrive at
\begin{align}
C_n^{(0)}-C_n^{(1)}&=(1-\varepsilon q^n+(\varepsilon q^n)^2)(C_{n-1}^{(0)}-C_{n-1}^{(1)})\notag\\[5pt]
&=\frac{1+\varepsilon q^{3n}}{1+\varepsilon q^n}(C_{n-1}^{(0)}-C_{n-1}^{(1)}).\label{c-15}
\end{align}
By \eqref{c-8} and \eqref{c-15}, we have
\begin{align}
C_n^{(0)}-C_n^{(2)}=\frac{1+\varepsilon q^{3n}}{1+\varepsilon q^n}(C_{n-1}^{(0)}-C_{n-1}^{(2)}).
\label{c-16}
\end{align}
By a repeated use of \eqref{c-15} and \eqref{c-16}, we obtain
\begin{align}
C_n^{(0)}-C_n^{(1)}=C_n^{(0)}-C_n^{(2)}=\frac{(-\varepsilon q^3;q^3 )_n}{(-\varepsilon q;q)_n}.
\label{c-17}
\end{align}

On the other hand, we have
\begin{align}
&\mathcal{L}_z\left((-\varepsilon zq;q)_n (-\varepsilon q/z;q)_n\right)
=\mathcal{L}_z\left(C_n^{(0)}+zC_n^{(1)}+z^2C_n^{(2)}\right)=C_n^{(1)}-C_n^{(2)},\label{c-18}\\[5pt]
&\mathcal{L}_z\left(z(-\varepsilon zq;q)_n (-\varepsilon q/z;q)_n\right)=
\mathcal{L}_z\left(zC_n^{(0)}+z^2C_n^{(1)}+C_n^{(2)}\right)=C_n^{(0)}-C_n^{(1)},\label{c-19}\\[5pt]
&\mathcal{L}_z\left(z^2(-\varepsilon zq;q)_n (-\varepsilon q/z;q)_n\right)=
\mathcal{L}_z\left(z^2C_n^{(0)}+C_n^{(1)}+zC_n^{(2)}\right)=C_n^{(2)}-C_n^{(0)}.\label{c-20}
\end{align}
Combining \eqref{c-8}, \eqref{c-17} and \eqref{c-18}--\eqref{c-20}, we complete the proof of \eqref{c-6-new}.

The proof of \eqref{c-7-new} runs analogously, and we omit the details.
\qed

\section{Proof of Theorems \ref{t-1} and \ref{t-2}}
{\noindent \it Proof of \eqref{a-1}.}
Letting $z\to \varepsilon z$ in \eqref{c-new-1}, we have
\begin{align*}
&\sum_{k=-n}^m (\varepsilon z)^{k+d} q^{\frac{k^2+sk}{2}}{n+m\brack n+k}\\[5pt]
&=(\varepsilon z)^{d-n}q^{\frac{n(n-s)}{2}}(1+\varepsilon zq^{-n+\frac{s+1}{2}})\cdots (1+\varepsilon zq^{-1})(1+\varepsilon z)(1+\varepsilon zq)\cdots (1+\varepsilon zq^{m+\frac{s-1}{2}})\\[5pt]
&=q^{\frac{1-s^2}{8}}(\varepsilon z)^{d-\frac{s+1}{2}}(1+\varepsilon z)(1+\varepsilon q^{n-\frac{s+1}{2}}/z)\cdots (1+\varepsilon q/z)(1+\varepsilon zq)\cdots (1+\varepsilon zq^{m+\frac{s-1}{2}}).
\end{align*}
Since $n-(s+1)/2\ge m+(s-1)/2$, we have
\begin{align}
&\sum_{k=-n}^m (\varepsilon z)^{k+d} q^{\frac{k^2+sk}{2}}{n+m\brack n+k}\notag\\[5pt]
&=q^{\frac{1-s^2}{8}}(\varepsilon z)^{d-\frac{s+1}{2}}(1+\varepsilon z)(1+\varepsilon q^{n-\frac{s+1}{2}}/z)\cdots (1+\varepsilon q^{m+\frac{s+1}{2}}/z)\notag\\[5pt]
&\times (1+\varepsilon q^{m+\frac{s-1}{2}}/z)\cdots (1+\varepsilon q/z)(1+\varepsilon zq)\cdots (1+\varepsilon zq^{m+\frac{s-1}{2}})\notag\\[5pt]
&=q^{\frac{1-s^2}{8}}(\varepsilon z)^{d-\frac{s+1}{2}}(1+\varepsilon z)(-\varepsilon q^{m+\frac{s+1}{2}}/z;q)_{n-m-s}\notag\\[5pt]
&\times (-\varepsilon zq;q)_{m+(s-1)/2}(-\varepsilon q/z;q)_{m+(s-1)/2}.\label{d-1}
\end{align}

For $r\in \{0,1,2\}$, let $C^{(r)}(q)$ be the sum of all the coefficients of $z^{3l+r} (l\in \mathbb{N})$ in
\begin{align*}
q^{\frac{1-s^2}{8}}(\varepsilon z)^{d-\frac{s+1}{2}}(1+\varepsilon z)(-\varepsilon q^{m+\frac{s+1}{2}}/z;q)_{n-m-s}.
\end{align*}
Taking the operator $\mathcal{L}_z$ to act on both sides of \eqref{d-1} and using \eqref{c-6}, we obtain
\begin{align}
&\sum_{k=-n}^m \varepsilon^{k+d}\left(\frac{k+d}{3}\right) q^{\frac{k^2+sk}{2}}{n+m\brack n+k}\notag\\[5pt]
&=\frac{(-\varepsilon q^3;q^3)_{m+(s-1)/2}}{(-\varepsilon q;q)_{m+(s-1)/2}}\left(C^{(1)}(q)-C^{(2)}(q)\right).\label{d-2}
\end{align}
Note that
\begin{align}
\mathcal{L}_z\left(q^{\frac{1-s^2}{8}}(\varepsilon z)^{d-\frac{s+1}{2}}(1+\varepsilon z)(-\varepsilon q^{m+\frac{s+1}{2}}/z;q)_{n-m-s}\right)=C^{(1)}(q)-C^{(2)}(q). \label{d-3}
\end{align}
Then the proof of \eqref{a-1} follows from \eqref{d-2} and \eqref{d-3}.
\qed

{\noindent \it Proof of \eqref{a-2}.}
Letting $z\to \varepsilon z$ in \eqref{c-new-1}, we have
\begin{align*}
&\sum_{k=-n}^m (\varepsilon z)^{k+d} q^{\frac{k^2+sk}{2}}{n+m\brack n+k}\\[5pt]
&=(\varepsilon z)^{d-n}q^{\frac{n(n-s)}{2}}(1+\varepsilon zq^{-n+\frac{s+1}{2}})\cdots (1+\varepsilon zq^{-\frac{1}{2}})(1+\varepsilon zq^{\frac{1}{2}})\cdots (1+\varepsilon zq^{m+\frac{s-1}{2}})\\[5pt]
&=q^{-\frac{s^2}{8}}(\varepsilon z)^{d-\frac{s}{2}}(1+\varepsilon q^{n-\frac{s+1}{2}}/z)\cdots (1+\varepsilon q^{\frac{1}{2}}/z)(1+\varepsilon zq^{\frac{1}{2}})\cdots (1+\varepsilon zq^{m+\frac{s-1}{2}}).
\end{align*}
Since $n-(s+1)/2\ge m+(s-1)/2$, we have
\begin{align}
&\sum_{k=-n}^m (\varepsilon z)^{k+d} q^{\frac{k^2+sk}{2}}{n+m\brack n+k}\notag\\[5pt]
&=q^{-\frac{s^2}{8}}(\varepsilon z)^{d-\frac{s}{2}}(1+\varepsilon q^{n-\frac{s+1}{2}}/z)\cdots (1+\varepsilon q^{m+\frac{s+1}{2}}/z)\notag\\[5pt]
&\times (1+\varepsilon q^{m+\frac{s-1}{2}}/z)\cdots (1+\varepsilon q^{\frac{1}{2}}/z)(1+\varepsilon zq^{\frac{1}{2}})\cdots (1+\varepsilon zq^{m+\frac{s-1}{2}})\notag\\[5pt]
&=q^{-\frac{s^2}{8}}(\varepsilon z)^{d-\frac{s}{2}}(-\varepsilon q^{m+\frac{s+1}{2}}/z;q)_{n-m-s}(-\varepsilon zq^{\frac{1}{2}};q)_{m+s/2}
(-\varepsilon q^{\frac{1}{2}}/z;q)_{m+s/2}.\label{d-4}
\end{align}
Using \eqref{c-7}, \eqref{d-4} and the same method as in the previous proof, we complete the proof of
\eqref{a-2}.
\qed

\vskip 5mm \noindent{\bf Acknowledgments.}
This work was supported by the National Natural Science Foundation of China (grant 12171370).


\begin{thebibliography}{99}

\small \setlength{\itemsep}{-.8mm}

\bibitem{andrews-b-1998}G.E. Andrews, The Theory of Partitions, Cambridge University Press, Cambridge, 1998.

\bibitem{berkovich-dm-2024}A. Berkovich, On the $q$-binomial identities involving the Legendre symbol modulo $3$, Discrete Math. 347 (2024), Art. 113824.

\bibitem{bu-jmaa-2022}A. Berkovich, A.K. Uncu, New infinite hierarchies of polynomial identities related to the Capparelli partition theorems, J. Math. Anal. Appl. 506 (2022), Art. 125678.

\end{thebibliography}
\end{document}